\documentclass[10pt,a4paper,oneside]{article}

\usepackage{latexsym}
\usepackage{amsbsy}
\usepackage{amsmath}
\usepackage{amssymb}
\usepackage{graphics}
\usepackage{epsfig}

\def\QQ{\mathbb{Q}}
\def\RR{\mathbb{R}}
\def\ZZ{\mathbb{Z}}
\def\CC{\mathbb{C}}
\def\i{\mathrm{int}}

\newtheorem{thm}{Theorem}[section]

\newtheorem{prop}[thm]{Proposition}
\newtheorem{cor}[thm]{Corollary}

\newenvironment{defin}
      {\vspace{2.5mm}\par\noindent
                \textbf{Definition.}}
             {\vspace{2.5mm}}

\newenvironment{proof}%
      {\par\noindent%
            \textbf{Proof.}}%
           { ~\hfill$\Box$\linebreak}

\newcommand{\dv}{ \! : \! }

\title{\large \bfseries SIMPLICIAL STRUCTURES OF KNOT COMPLEMENTS}

\author{\normalsize ALEKSANDAR MIJATOVI\'{C}}

\date{}

\renewenvironment{abstract}
            {\begin{quotation}\noindent\small
                              \textsc{Abstract}.\hspace{0.5mm}}
            {\end{quotation}}

\begin{document}

\maketitle

\begin{abstract}
It was shown in~\cite{mijatov2}
that there exists an explicit bound 
for the number of Pachner moves needed to connect any two triangulation of
any 
Haken 3-manifold which contains no fibred sub-manifolds as 
strongly simple pieces of its JSJ-decomposition.
In this paper we prove a generalisation
of that result to all knot complements.
The explicit formula for the bound is in terms of the numbers of 
tetrahedra in the two triangulations. 
This gives a conceptually
trivial algorithm for recognising any knot complement
among all 3-manifolds.
\end{abstract}

\begin{center}
\section{\normalsize \scshape INTRODUCTION}
\label{sec:intro}
\end{center}

The main aim of this paper is to understand the relation between
any two triangulations of a given
knot complement.
In order to alter one simplicial structure 
(throughout the paper the term
\textit{simplicial structure} 
will be used as a synonym for a triangulation) to another 
we need the 
following definition.

\begin{figure}[!hbt]
  \begin{center}
    \epsfig{file=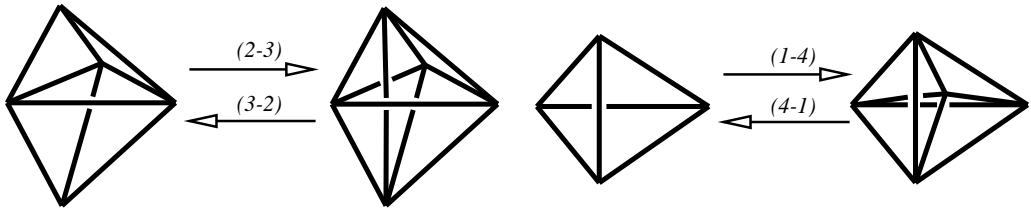}
        \caption{{\small Three dimensional Pachner moves.}}
        \label{fig:3pm}
   \end{center}
\end{figure}

\begin{defin}
Let
$T$
be a triangulation of a compact PL
\mbox{$n$-manifold}
$M$.
Suppose
$D$
is a combinatorial
\mbox{$n$-disc}
which is a sub-complex both of
$T$
and of the boundary of a standard
\mbox{$(n+1)$-simplex}
$\Delta^{n+1}.$
A
\textit{Pachner move}
consists of changing
$T$
by removing the sub-complex
$D$
and inserting
$\partial\Delta^{n+1}-\i(D)$
(for
$n$
equals 3, see figure~\ref{fig:3pm}).
\end{defin}

It is an immediate consequence of this definition that there
are precisely 
$(n+1)$
possible Pachner moves in dimension
$n$.
If our 
$n$-manifold 
$M$
has non-empty boundary, then the moves 
from this definition do not alter the induced triangulation of
$\partial M$.
But changing the simplicial structure 
of the boundary with an
$(n-1)$-dimensional Pachner move can be achieved by
gluing onto (or removing from) our manifold
$M$
the standard
$n$-simplex
$\Delta^n$
that exists by the definition of the move.

The moves from the above 
definition are in some sense completely general.
It was proved by Pachner in~\cite{pachner} that any two
triangulations of the same PL
$n$-manifold are related by a finite sequence of Pachner
moves and simplicial isomorphisms. It is well known
(see proposition 1.3 in~\cite{mijatov})
that in case of a fixed 3-manifold
$M$
a computable function, depending only on the
number of tetrahedra in the triangulations
of
$M$, 
bounding the
length of the sequence from Pachner's theorem, gives
an algorithm for recognising
$M$
among all 3-manifolds. The following theorem
gives an explicit formula for such a bound in case
$M$
is a 
knot complement.\\

\begin{thm}
\label{thm:knot}
Let
$P$
and
$Q$
be two triangulations of a knot complement
$M$
which contain
$p$
and
$q$
tetrahedra respectively. Then there exists a sequence of
Pachner moves of length at most
$e^{2^{ap}}(p)+ e^{2^{aq}}(q)$
which transforms
$P$
into a triangulation isomorphic to
$Q$.
The constant
$a$
is bounded above by
$200$.
The homeomorphism of
$M$
that realizes this simplicial isomorphism is supported in
the characteristic sub-manifold
of
$M$.
\end{thm}

The triangulations appearing in theorem~\ref{thm:knot} are 
allowed to be non-combinatorial, which means that the simplices
are not (necessarily) uniquely determined by their vertices.       
The exponent in the above expression containing the exponential
function
$e(x)=2^x$
stands for the composition of the function with itself rather 
than for multiplication. Since the formula in theorem~\ref{thm:knot}
is explicit, it gives a conceptually trivial algorithm 
to recognise
any knot complement 
among all 3-manifolds
(just make all possible
sequences of Pachner moves whose length is smaller than the bound!). 

In addition to the recognition problem for all knot complements,
theorem~\ref{thm:knot} gives
a simple
procedure
that
can be used to decide if any knot, represented by a knot
diagram, is the same as our given knot. This procedure
will be described in section~\ref{sec:main}.

Theorem~\ref{thm:knot} can be obtained as a direct consequence
of theorem~\ref{thm:final} which is contained in
the next section. The proof uses
the main result from~\cite{mijatov2}. This however is not enough
because we also need to deal with fibred knot complements. It
has been know for a long time that an algorithm capable of
deciding whether two homeomorphisms of the fibre are 
conjugate can be used to recognise fibred manifolds. 
Our approach, however, is completely different.
We avoid this problem by applying the crucial 
proposition~\ref{prop:knot} which ensures that the 
first surface in the canonical hierarchy is not a
fibre. In other words this means that the procedure
used to prove theorem 3.1 in~\cite{mijatov2} also 
works in this setting and implies theorem~\ref{thm:knot}.
The proof of proposition~\ref{prop:knot} depends heavily
on the existence of a separating incompressible surface
which was established by Culler and Shalen in~\cite{culler}.

In section~\ref{sec:main} 
we give a precise description
of the class of 3-manifolds
which appear in theorem~\ref{thm:final}. After stating the main
theorem we also outline its proof.
Section~\ref{sec:separating}
contains the relevant results from~\cite{culler}
and their application to our setting. 
In the last section we prove the key propositions~\ref{prop:2bc}
and~\ref{prop:knot}.

\begin{center}
\end{center}
\begin{center}
\section{\normalsize \scshape THE MAIN THEOREM}
\label{sec:main}
\end{center}

Let's start by recalling some well-known definitions. A
3-manifold 
$M$
is \textit{irreducible} if every embedded 2-sphere
in it bounds a 3-ball. 
A properly embedded surface
$F$
in an irreducible 3-manifold
$M$
is
\textit{injective} if the homomorphism
$\pi_1(F)\rightarrow \pi_1(M)$,
induced by the inclusion of
$F$
into
$M$,
is a monomorphism. A surface
$F$
is said to be \textit{incompressible} if no component of
it
is a 2-sphere 
or a disc
and if for every disc
$D$
in
$M$
with
$D\cap S=\partial D$,
there is a disc
$D'$
in
$S$
with $\partial D=\partial D'$.

There is also a relative notion of incompressibility which we
will need to consider. A surface
$F$
is
$\partial$-\textit{incompressible} if for each disc
$D$
in
$M$,
such that
$\partial D$
splits into two arcs
$\alpha$
and
$\beta$
meeting only at their common endpoints
with
$D\cap F=\alpha$
and
$D\cap \partial M=\beta$,
there is a disc
$D'$
in
$F$
with
$\alpha\subset\partial D'$
and
$\partial D'-\alpha\subset \partial F$.

A \textit{horizontal boundary} of an
$I$-bundle over a surface is the part of the boundary
corresponding to the
$\partial I$-bundle. The \textit{vertical boundary} is
the complement of the horizontal boundary and consists of annuli
that fibre over the bounding circles of the base surface.
It is a well-known fact that a properly embedded one-sided surface
in
$M$
is injective if and only if the horizontal boundary of its regular
neighbourhood is incompressible.
An 
irreducible 3-manifold 
$M$
with possibly empty
incompressible boundary 
is
\textit{Haken} 
if it contains an injective
surface different from a disc or a 2-sphere. 
A torus (resp. annulus) in 
$M$
that is incompressible and is not boundary parallel is sometimes
referred to as an \textit{essential torus} (resp.
\textit{essential annulus}).
Notice that in a Haken 
3-manifold an essential annulus can not be
$\partial$-compressible.

Before stating theorem~\ref{thm:final}
we need to recall some more standard terminology.
A \textit{surface bundle} with an orientable fibre
$S$
is just a mapping torus, i.e. a quotient
$S\times I/(x,0)\sim(\varphi(x),1)$,
for some orientation preserving surface automorphism
$\varphi\dv S\rightarrow S$.
Since
$S$
is orientable
this construction gives an orientable 3-manifold. But for a non-orientable
surface
$R$
with a non-trivial two-sheeted covering
$S\rightarrow R$,
the mapping cylinder of the covering projection is an
orientable twisted
$I$-bundle over
$R$.
Gluing two such
$I$-bundles
together along their horizontal boundaries
by an automorphism of
$S$
gives a 3-manifold
$N$
which is foliated by parallel copies of
$S$
and the two copies of
$R$.
The leaves of this foliation are the
``fibres'' of a natural projection map
$N\rightarrow I$,
where the two copies of
$R$
are the pre-images of the endpoints of the interval
$I$.
Such a 3-manifold
$N$
will be called a \textit{semi-bundle} (with fibre
$S$) over an interval
$I$.
The surfaces
$S$
and
$R$
can be either closed or bounded.

Manifolds which are homeomorphic to semi-bundles
do sometimes arise naturally.
The simplest
example is a connected sum of two projective spaces
$\RR P^3\#\RR P^3$
where the fibre is a 2-sphere.
On the other hand a semi-bundle structure can never arise in a knot
complement. This is because the boundary circles of the
two non-orientable leaves would be disjoint curves in the
boundary torus and
could therefore be capped off by the
annuli they bound in the torus. This would then
give a closed non-orientable
surface in
$S^3$.
This observation will be crucial for us because it will
insure that theorem~\ref{thm:final} implies theorem~\ref{thm:knot}.

The JSJ-decomposition of 
a Haken 3-manifold consists
of strongly simple pieces,
$I$-bundles
and
Seifert fibred spaces
(for precise definitions see
section 2 in~\cite{mijatov2}). The strongly simple pieces are the ones 
that contain all the interesting topological information
about 
$M$
but also have the crucial property 
of being both atoroidal and an-annular (i.e.
all incompressible annuli and
tori in them are boundary parallel). Loosely speaking the union
of all the components of the JSJ-decomposition
that are either homeomorphic to 
$I$-bundles or to Seifert fibred spaces constitute the characteristic
sub-manifold 
$\Sigma$
of 
$M$. 
Before we state theorem~\ref{thm:final}, we should remind ourselves that
the exponent in the formula below, containing the exponential function
$e(x)=2^x$,
stands for the composition of the function with itself rather
than for multiplication.\\

\begin{thm}
\label{thm:final}
Let
$M$
be a Haken 3-manifold. Assume that the
strongly simple pieces of its JSJ-decomposition
are not homeomorphic to any of the following
types of 3-manifolds:
a closed surface
semi-bundle which is a rational homology 3-sphere,
a closed surface bundle with the first Betti number equal to
one or a surface bundle with a single boundary component which
contains no closed injective surfaces (other than
the boundary torus) and which
is at the same time homeomorphic
to a surface semi-bundle.
Let
$P$
and
$Q$
be two triangulations of
$M$
that contain
$p$
and
$q$
tetrahedra respectively. Then there exists a sequence of
Pachner moves of length at most
$e^{2^{ap}}(p)+ e^{2^{aq}}(q)$
which transforms
$P$
into a triangulation isomorphic to
$Q$.
The constant
$a$
is bounded above by
$200$.
The homeomorphism of
$M$,
that realizes this simplicial isomorphism, is supported in
the characteristic sub-manifold
$\Sigma$
of
$M$
and it does not permute the components of
$\partial M$.
\end{thm}

Since knot complements satisfy the hypothesis of
theorem~\ref{thm:final} we get a conceptually
trivial algorithm for determining whether
any 3-manifold is homeomorphic to a complement 
of a given knot in 
$S^3$.
Moreover theorem~\ref{thm:final} also gives a simple 
procedure to determine whether
any knot diagram 
represents a knot 
which is isotopic to our given knot. It is enough to
establish whether their respective complements are homeomorphic 
(which we already know how to do) and,
if they are, to determine whether the homeomorphism maps the meridian
of one onto the meridian of the other.
If the boundary torus of the knot complement is
not contained in the characteristic sub-manifold, then
the homeomorphism
from theorem~\ref{thm:final} equals
the identity
on the boundary. If on the other hand the bounding torus is contained in
$\Sigma$,
then we first make sure that the simplicial structures on the boundary of
both knot complements coincide.
It follows from the proof of theorem 3.1 in~\cite{mijatov2} 
that this is enough to make our homeomorphism equal   
to the identity on the boundary torus.
So in this way, using theorem~\ref{thm:final}, we can solve the
recognition problem for any knot.

The proof of theorem~\ref{thm:final} follows the same
lines as the proof of the main theorem in~\cite{mijatov2}.
The main tool for probing the topology of the simple pieces
of the JSJ-decomposition of 
$M$
is the canonical hierarchy (see section 4 in~\cite{mijatov2}).
This hierarchy is based on Haken's original recognition
algorithm for non-fibred 3-manifolds which are sufficiently
large. 
The last step of Haken's program for the
classification of
sufficiently large 3-manifolds is the solution of 
the recognition problem for
surface bundles. Haken knew that an algorithm
capable of deciding whether two automorphisms
of the fibre are conjugate would suffice. The
algorithmic solution of the conjugacy
problem in the mapping class group of the fibre was first
proved by Hemion in~\cite{hemion1} and has been reproved
many time since. We, however, take a completely different
approach in the fibred case situation. 

The proof of theorem 3.1 in~\cite{mijatov2} starts by 
constructing the canonical hierarchy (section 4 of~\cite{mijatov2})
in 
$M$.
The first surface 
$S_1$
in the hierarchy consists of the JSJ-system (i.e. canonical tori
and annuli, see section 2 in~\cite{mijatov2}) and of the 
closed two-sided 
injective surfaces in the 
strongly simple pieces of the
JSJ-decomposition of
$M$. 
The surface 
$S_1$
is defined so that the complement 
$M-\mathrm{int}(\mathcal{N}(S_1))$
contains no closed orientable incompressible surfaces which
are not boundary parallel. 
In each component of 
$M-\mathrm{int}(\mathcal{N}(S_1))$,
which is disjoint from the characteristic sub-manifold
$\Sigma$,
we take the surface 
$S_2$
to be
a bounded two-sided incompressible surface with
the largest Euler characteristic in that piece.
It was shown in section 4 of~\cite{mijatov2} that the
components of
$M-\mathrm{int}(\mathcal{N}(S_1\cup S_2))$
are topologically equivalent to
$I$-bundles, handlebodies and compression bodies. We continue by
cutting these complementary regions using step 3 of the
canonical hierarchy. 
The key lemma 4.2 of~\cite{mijatov2} tells us that the 
canonical hierarchy decomposes 
$M$
in a manageable way if and only if no component of
$S_1$
or 
$S_2$
is a fibre in a bundle structure or a semi-bundle structure of a 
simple piece in 
the JSJ-decomposition of
$M$. 
In~\cite{mijatov2} we made sure that this was not 
the case by hypothesising away all 3-manifolds that
contain strongly simple pieces which support bundle
and semi-bundle structures.

In theorem~\ref{thm:final} we allow for many of the 3-manifolds
from the ``fibred''  family. We avoid problems by making sure that the
crucial components of 
$S_1$
and
$S_2$
are not fibres. 
Once we show that such surfaces exist and that they
have bounded
normal complexity, 
everything works in exactly the same way 
as in the proof of theorem 3.1 in~\cite{mijatov2}. 
If a fibred strongly simple piece of 
$M$
has boundary,
we use 
propositions~\ref{prop:2bc} and~\ref{prop:knot}.
Their proofs 
are not based on the solution
of the conjugacy problem. Instead they use a different deep
(geometric) fact from~\cite{culler} which says that our
surface bundle has to contain a separating incompressible
surface. 
The same philosophy of looking for a surface
which is not a fibre can be applied to 
all closed atoroidal
3-manifolds that have enough homology. The existence of
such a surface is guaranteed by the work of Thurston 
in~\cite{thurston}. A bound on the normal complexity
of such a surface can be obtained directly from
the results of Wang and Tollefson 
in~\cite{wang}. A more detailed description of this procedure
will be given in the last section.

Once we find the surfaces that are not fibres, 
we apply the canonical hierarchy techniques
(section 4 in~\cite{mijatov2})
together with theorem 1.2 of~\cite{mijatov} to
all strongly simple pieces
of
$M$.
This makes it
possible to connect any two triangulations of the simple sub-manifolds
by a sequence of Pachner moves.
The subdivision of the original triangulation in the
characteristic sub-manifold can be altered directly by applying
the main theorem
of~\cite{mijatov1}.

The reason why our strategy fails for some surface bundles
with a single boundary component is the following.
There seems to be no way of
ensuring that the component of the surface
$S_2$
(we are trying to construct) 
is neither a fibre in the bundle structure nor in
the semi-bundle structure of the piece. Since 3-manifolds
supporting both of these fibred structures exist, we have to
exclude them by hypothesis.

It seems that dealing with triangulations of 
closed fibred manifolds which do not satisfy the assumptions of 
theorem~\ref{thm:final} requires solving
the conjugacy problem in the mapping class
group of the fibre. On the other hand theorem~\ref{thm:final}
can be used to solve the conjugacy problem 
for the elements in
the mapping class
group of a surface with at least two punctures,
which fix the 
boundary circles. This is because
any two orientation preserving homeomorphisms of a
surface are conjugate if and only if the two associated mapping
tori are homeomorphic via a homeomorphism
which maps fibres to fibres. If the surface bundle has at least two
boundary components we can check, using theorem~\ref{thm:final},
if the mapping tori are homeomorphic. While we are changing
one of the triangulations using Pachner moves, we can at the same time
keep track of the original fibre in normal form. This is
because at the beginning the fibre can be easily isotoped 
into normal form with respect to the starting triangulation.
It follows directly from the definition of Pachner moves
that we can keep it in normal form after we make each move.
If in the end the two surface bundles are homeomorphic, we must
check whether the fibres we have been keeping track of are isotopic.
Since this can be done algorithmically,
because both surfaces are represented by their normal forms, 
we can use ti to solve the
conjugacy problem.

\begin{center}
\end{center}
\begin{center}
\section{\normalsize \scshape SEPARATING INCOMPRESSIBLE SURFACES}
\label{sec:separating}
\end{center}

The following amazing result is one of the main theorems of~\cite{culler}.
We will use it to prove the existence of surfaces which are not fibres
in bounded 3-manifolds.\\

\begin{thm}[Culler, Shalen]
\label{thm:cs}
Let 
$M$
be a compact connected orientable 3-manifold whose boundary is
a non-empty union of tori. Suppose that 
$H_1(\partial M;\QQ)\rightarrow H_1(M;\QQ)$
is surjective and that 
$M$
is not homeomorphic to a solid torus or
$(\mathrm{torus})\times I$.
Then for each component 
$B$
of
$\partial M$,
there is a separating connected incompressible surface in
$M$,
which is not 
$\partial$-parallel and whose boundary is not empty and
is contained in
$B$.
\end{thm}

This is a deep fact indeed. The starting point of its proof is
Thurston's geometrization of simple bounded 3-manifolds. It then 
applies some algebro-geometric techniques to analyse a certain 
complex curve in the set of characters of representations of
$\pi_1(M)$
in 
$SL_2(\CC)$.
An ``ideal point'' point on this curve gives rise to a 
``splitting'' of the fundamental group which in turn 
can be used to produce an incompressible surface in
$M$
that is not boundary parallel. We will now apply theorem~\ref{thm:cs}
to our setting.
\vspace{1.5mm}

\begin{cor}
\label{cor:cs}
Let 
$M$
be a connected irreducible atoroidal 3-manifold with 
boundary that is a non-empty collection of tori. 
Assume also that
$M$
is not a Seifert 
fibred space and that it does not contain a closed
injective surface which is not
boundary parallel. Then the following holds.
\newcounter{h}
\begin{list}{(\alph{h})}{\usecounter{h}\parsep 0mm\topsep 1mm}
\item Assume that 
      $\partial M$
      is disconnected and that
      $B$
      is one of its components. 
      Then 
      $M$
      contains
      a connected orientable separating incompressible
      $\partial$-incompressible surface 
      with non-empty boundary which is neither a
      fibre in a bundle structure over a circle nor in a semi-bundle
      structure over an interval. Also the boundary of this surface
      is contained in
      $B$.     

\item Assume that
      $\partial M$
      is a single torus and that 
      $M$
      is not homeomorphic to a surface semi-bundle.
      Then 
      $M$
      contains 
      a connected orientable separating incompressible
      $\partial$-incompressible surface with non-empty 
      boundary which is not a
      fibre in any bundle structure over a circle.
\end{list}
\end{cor}

This corollary can be applied directly to any strongly 
simple piece in the JSJ-decomposition of our manifold
which contains no components of the surface 
$S_1$
after step 1 of the canonical hierarchy and
whose boundary consist of tori. 
We have seen in previous section
that knot
complements in 
$S^3$
do not admit a semi-bundle structure. So corollary~\ref{cor:cs}
applies to all knot complements. 
\vspace{2mm}

\begin{proof}
It is well-known that every non-trivial element of
$H_2(M;\QQ)$
gives rise to a closed orientable
non-separating surface in
$M$.
Since such surfaces do not exist in
$M$
this
implies that
$H_2(M;\QQ)$
must be trivial. By Poincar\'{e} duality it follows that 
$H_1(M,\partial M,;\QQ)$
is also trivial. The exact homology sequence of the pair
$(M, \partial M)$
gives that the homomorphism 
$H_1(\partial M;\QQ)\rightarrow H_1(M;\QQ)$
is surjective.
Now we can apply theorem~\ref{thm:cs} to our manifold
$M$.

If
$\partial M$
is disconnected, we obtain a separating orientable 
connected incompressible
surface 
$F$
in
$M$
which is not boundary parallel and which is disjoint from 
at least one component of
$\partial M$.
Also 
$\partial F$
is 
not empty and is
contained in
$B$.
The surface 
$F$
can not be a fibre because it is disjoint from
$(\partial M)-B$.
We need to show that it is
$\partial$-incompressible. 
Since 
$F$
is separating it contains at least two boundary circles
in 
$B$.
The bounding circle of a 
$\partial$-compression disc for 
$F$
is a union of two arcs: one in 
$F$
and one in
$B$.
The arc in
$B$
either runs between two distinct components of
$\partial F$
or hits a single circle in
$\partial F$
twice from the same side.
In the former case
we can construct, using the annulus in
$B$
between 
the two circles of
$B\cap F$,
a genuine compression disc for
$F$. 
Since 
$F$
is incompressible and 
$M$
is irreducible this would imply that 
$F$
is
$\partial$-parallel. 
In the latter case it is even easier to construct a compression
disc for
$F$
which
again leads into contradiction.
So 
$F$
must be 
$\partial$-incompressible. 

If 
$\partial M$
consists of a single torus, then theorem~\ref{thm:cs}
gives us a surface with the same properties as 
$F$,
which can not be a fibre in a bundle structure over
a circle because it is separating.
\end{proof}

\begin{center}
\end{center}
\begin{center}
\section{\normalsize \scshape NON-FIBRES IN NORMAL FORM}
\label{sec:non-fibre}
\end{center}

In this section 
we are going to describe how to construct an incompressible 
surface, 
which is
not a fibre, 
in a triangulated (semi-)bundle 
$M$
that satisfies the assumptions of corollary~\ref{cor:cs}.
There are two cases depending on whether 
$\partial M$
is connected or not. We will first deal with the latter case
which is easier. At the end of this section we will discuss the closed
case and thus complete the proof of theorem~\ref{thm:final}.
Let's start by recalling some normal surface theory.

Let
$T$
be a triangulation of a 3-manifold
$M$.
An arc in a 2-simplex of
$T$
is
\textit{normal}
if its ends lie in different
sides of a 2-simplex. A simple closed
curve in the 2-skeleton of
$T$
is a \textit{normal curve} if it intersects each 2-simplex of
$T$
in normal arcs. A properly embedded surface
$F$
in
$M$
is in \textit{normal form} with respect to
$T$
if it intersects each tetrahedron in
$T$
in a collection of discs all of whose boundaries are normal
curves consisting of 3 or 4 normal arcs, i.e. triangles and
quadrilaterals.
A \textit{normal disc} is a triangle or a quadrilateral.
There are precisely seven normal disc types in any tetrahedron
of
$T$.
An isotopy of
$M$
is called a \textit{normal isotopy} with respect to
$T$
if it leaves all simplices of
$T$
invariant. In particular this implies that it is fixed on the vertices of
$T$.

A normal surface is determined, up to normal isotopy, by the number
of normal disc types in which it meets the tetrahedra of
$T$.
It therefore defines a vector with
$7t$
coordinates.
Each coordinate represents the number of copies of  normal disc
types that are contained in the surface
($t$
is the number of tetrahedra in
$T$).
It turns out
that there is a certain restricted linear system that such
a vector is a solution
of. Moreover there is a one to one
correspondence between
the solutions of that restricted linear system and normal
surface
in
$M$.
If the sum of two vector solutions of this
system satisfies the restrictions on the system, then it represents a normal
surface in
$M$.
On the other hand there is a geometric process called
\textit{regular alteration} (see figure 2 in~\cite{bart})
which can be carried out
on the normal surfaces representing the summands and which
yields the normal surface corresponding to the sum.
It follows directly from the definition of regular alteration
that the Euler characteristic is additive over normal addition.
We can define the
\textit{weight}
$w(F)$
of a surface
$F$,
which is transverse to the 1-skeleton of
$T$,
to be the number of points of intersection between
the surface and the 1-skeleton. Since regular alteration only changes
the surfaces involved away from the 1-skeleton,
the weight too is additive over normal addition.

A normal surface is called
\textit{fundamental} if the vector corresponding to
it is not a sum of two integral solutions of the linear
system.
The solution space of the linear system projects down to a
compact convex linear cell which is called the
\textit{projective solution space}.
A \textit{vertex surface}
is a connected two-sided normal surface that projects onto a vertex
of the projective solution space (see~\cite{wang} for a more
detailed description).
The normal sum
$F=F_1+F_2$
is in \textit{reduced form} if the number of components of
$F_1\cap F_2$
is minimal among all normal surfaces
$F_1'$
and
$F_2'$
isotopic 
to 
$F_1$
and
$F_2$
respectively
such that
$F=F_1'+F_2'$.

Before we proceed we need to define
two (very simple) kinds of complexities of the surfaces 
embedded in our triangulated 3-manifold
$M$.
First there is the \textit{normal complexity},
i.e. the number of normal pieces a minimal
weight representative in the isotopy class of the surfaces consist of.
Second
there is the \textit{topological complexity} of a surface 
which is defined in terms of its components in the following way.
To each component we assign its negative Euler characteristic and
then define the complexity
to be the sum over all of the components.
Since there are
no 2-spheres, discs or projective planes among the surfaces
we are trying to construct in this paper, 
their topological complexity will coincide with
the Thurston complexity as defined in~\cite{thurston}.
Now we can state the following proposition.\\

\begin{prop}
\label{prop:2bc}
Let 
$M$
be a triangulated 3-manifold which satisfies the assumptions
of corollary~\ref{cor:cs} (a) and which is also
an-annular. 
Let
$B$
be a component of
$\partial M$. 
Let
$F$
be a two-sided connected incompressible 
$\partial$-incompressible surface in 
$M$
with
non-empty boundary lying in
$B$,
which
minimises the topological complexity among all
such surfaces. Then we can isotope 
$F$
into normal form so that it is a sum of at most two
fundamental surfaces. 
\end{prop}

\begin{proof}
Corollary~\ref{cor:cs} 
guarantees the existence of 
at least one surface with the properties from the proposition.
That surface is also separating while our
$F$
might not be. Assume that 
$F$
is in normal form and that it has minimal weight in its
isotopy class. Now
express 
$F$
as a sum of fundamental surfaces: 
$F=k_1F_1+\ldots+k_nF_n$.
By theorem 2.3 from~\cite{mijatov1} we can conclude that each
$F_i$
is incompressible and
$\partial$-incompressible. In fact all the summands are
injective because we can apply the same theorem
to 
$2F$
(which also minimises the weight in its isotopy class since 
$F$
is two-sided). 

If some
$F_i$
were closed, then it would have to be a
$\partial$-parallel torus
by our assumption on 
$M$.
We would then get
$F=F_i+S$
where 
$S$
is some normal surface in
$M$.
We can assume that the 
sum
$F=F_i+S$
is in reduced form. 
Lemma 2.2 from~\cite{mijatov1} then 
implies that no component of the surface
$F-(F_i\cap S)$
is a disc. 
So the space
$F_i\cap S$
is a 1-manifold that is homeomorphic to
a disjoint union of non-trivial parallel
simple closed curves in the torus
$F_i$.
Let
$X$
be the
$(\mathrm{torus})\times I$
region between
$F_i$
and the toral boundary component of
$M$.
Then the components of the surface
$S\cap X$
must be injective in
$X$
simply because the patches of
$F=F_i+S$
are injective by lemma 2.2 from~\cite{mijatov1} and
$\partial X$
is incompressible in
$M$.
So, since
$S\cap X$
contains no closed components, it
consists of incompressible annuli that are either
disjoint from
the torus
$X\cap \partial M$
or are spanning annuli in
$X$.
Let
$B$
be an outermost annular component of
$S\cap X$
which is disjoint from
$X\cap \partial M$
and let
$A$
be the annulus in
$F_i$
that is parallel to
$B$.
There are three possible (essentially different)
ways  normal alteration can
act on
$\partial A$. They are depicted by figure~\ref{fig:alter}.

\begin{figure}[!hbt]
\vspace{-2cm}
 \begin{center}
   \vspace{2cm}
    \epsfig{file=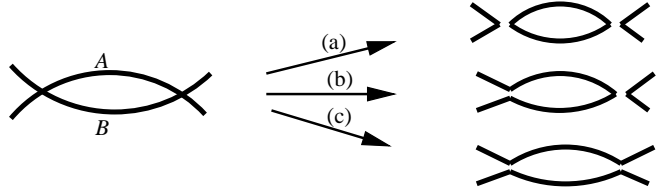}
  \caption{\small Possible normal alterations.}
  \label{fig:alter}
 \end{center}
\end{figure}

They all lead to contradiction. Case (a) produces a disconnected sum.
In case (b) we can isotope the union of the patches
$A$
and
$B$
over the solid torus that they bound, to reduce the weight of
$F$.
If both
$A$
and
$B$
had zero weight, then there would exist a normal isotopy that would
reduce the number of components in
$F_i\cap S$.
This contradicts the reduced form assumption. Case (c) contradicts it
as well, because the surfaces we obtain after we do the normal alterations
along
$\partial A$,
are isotopic to
$F_i$
and
$S$,
but have fewer components of intersection. 

So the only possible components of
$S\cap X$
are the spanning annuli,
i.e. the ones that are vertical in the product
structure of
$X$.
There are essentially only two different ways
of doing normal alterations along all the simple closed curves in
$F_i\cap S$
if we want to obtain a connected surface.
They lead to contradiction because the surface 
$F$
we get in both cases is isotopic to 
$S$.
This contradicts the assumption that 
$F$
has
minimal weight in its isotopy class. 

We have just shown that 
$\partial F_i$
is not empty and that it is contained in
$B$
for every 
$i=1,\ldots,n$.
Since 
$M$
is an-annular, no
$F_i$
can be an annulus or a Moebius band. Since all
$F_i$
are connected, none of them can be two-sided (unless the sum
has only one summand). But if we have at least three
one-sided surfaces in the sum, then the double of one of them
is going to satisfy all the necessary conditions and will 
have its topological complexity smaller than that of
$F$.
Hence the proposition follows.
\end{proof}

Notice that the assumption that 
$M$
is an-annular does not create any problems for
us because the JSJ-pieces we are interested in have this
property by lemma 4.2 in~\cite{mijatov2}. We are now going to deal
with the remaining case when our
manifold has only one boundary component.

\begin{prop}
\label{prop:knot}
Let 
$M$
be a triangulated 3-manifold which satisfies the assumptions
of corollary~\ref{cor:cs} (b) and which is also
an-annular. Let 
$t$
be the number of tetrahedra in 
$M$. 
Assume further that 
$M$
is a surface bundle over a circle. Let 
$F$
be a connected two-sided separating incompressible 
$\partial$-incompressible surface with non-empty 
boundary that has the smallest
topological complexity among all such surfaces. 
Then 
$F$
can be isotoped into normal form so that it consists
of not more than
$2^{40t}$
normal discs.
\end{prop}

It is clear that the surface 
$F$
from proposition~\ref{prop:knot}
can not be a fibre in the bundel structure of
$M$
because it is separating. Also by applying proposition~\ref{prop:knot}
while constructing components of the surface
$S_2$
in step 2 of the canonical hierarchy, we can make sure we never
adjoin a fibre. So if our 3-manifold has a non-trivial 
JSJ-decomposition, the only way we can run into trouble is if
there is a strongly simple piece with a single toral
boundary component
which is a semi-bundle 
and a bundle at the same time.
So we can use
propositions~\ref{prop:2bc} and~\ref{prop:knot},
combined with the proof of theorem 3.1 in~\cite{mijatov2}, 
to show that theorem~\ref{thm:final} holds. 
\vspace{2mm}

\begin{proof}
Like in the proof of proposition~\ref{cor:cs}
we have
$H_1(M,\partial M;\QQ)\cong H_2(M;\QQ)=0$. 
The exact homology sequence of the pair 
$(M,\partial M)$
yields a short exact sequence
$$0\rightarrow H_2(M,\partial M;\QQ)\rightarrow H_1(\partial M;\QQ)
\rightarrow H_1(M;\QQ) \rightarrow 0$$
of vector spaces over the field
$\QQ$.
Since
$H_2(M,\partial M;\QQ)\cong H_1(M;\QQ)$ 
by Poincar\'{e} duality
and 
$\mathrm{dim}(H_1(\partial M;\QQ))=2$, 
we get that 
$\beta_1(M;\QQ)=1$.
Again 
Poincar\'{e} duality implies that
$H_2(M,\partial M;\ZZ)$
is isomorphic to
$H^1(M;\ZZ)$ 
and is therefore torsion-free. In other words we have
shown that
$H_2(M,\partial M;\ZZ)$
is isomorphic to 
$\ZZ$,
which will be very useful later on.

Let's now isotope 
$F$
into normal form so that it minimises the weight
in its isotopy class. We can now express it as 
a sum of fundamental surfaces in the usual way:
$F=k_1F_1+\ldots+k_nF_n$.
Like in the proof of proposition~\ref{prop:2bc}
we can conclude that each 
$F_i$
is an injective 
$\partial$-incompressible bounded surface
with 
$\chi(F_i)<0$.
The Euler characteristic inequality comes 
from the fact that 
$M$
is an-annular and can therefore not contain non-trivial
annuli and Moebius bands. The same argument shows that
any normal surface which appears as a summand 
of
$F$
has to be an injective
$\partial$-incompressisble surface with non-empty
boundary and strictly negative Euler characteristic.
Clearly no
$F_i$
can be two-sided and separating
(unless 
$F$
itself is fundamental). But if
a surface
$F_j$
is one-sided then its double 
$2F_j$
is a bounded connected 
two-sided separating incompressible 
$\partial$-incompressible surface with negative 
Euler characteristic that satisfies the inequality
$$2\chi(F_j)=\chi(2F_j)\leq\chi(F)=k_1\chi(F_1)+\ldots+k_n\chi(F_n)\leq
 -\sum_{i=1}^nk_i.$$
The first inequality comes from the fact that
$F$
minimises topological complexity in the family of surfaces
which contains
$2F_j$
and the second one is simply saying that 
no 
$F_i$
is an annulus or a Moebius band.
The well-known bound on the number of triangles
and quadrilaterals in a fundamental normal surface
(see lemma 6.1 in~\cite{hass}) implies the inequality
$-2\chi(F_j)<2^{20t}$.
So it follows that
the sum representing 
$F$
has at most
$2^{20t}$
fundamental summands. This proves the proposition 
in case one of the surfaces 
$F_i$
is 
one-sided.
The following claim is central in all that follows. 
\vspace{2mm}

\noindent{\textbf{Claim.}} Let the surface
$S$
be a fibre of the bundle 
$M$.
Then any connected two-sided incompressible 
$\partial$-incompressible surface 
$R$
in 
$M$,
which is not separating,
is actually isotopic to
$S$.
\vspace{2mm}

It is clear that 
$R$
represents a non-trivial element in 
$H_2(M,\partial M;\ZZ)=\ZZ$.
This element is primitive because 
$R$
is connected (see lemma 1 in~\cite{thurston}). 
So if we choose our orientations correctly, we get that 
the surfaces 
$R$
and
$S$
represent the same class in
$H_2(M,\partial M;\ZZ)$.
Now
we want to lift 
$R$
to the infinite cyclic cover of
$M$
which corresponds to the element
$[S]$
in 
$H^1(M;\ZZ)$.
This covering space is clearly homeomorphic to
$S\times \RR$.
The surface 
$R$
lifts if and only if every element in
$\pi_1(R)$,
when viewed as a homology class,
has trivial algebraic intersection with
$[S]$.
This condition is satisfied because
the surfaces 
$R$
and
$S$
are homologous and 
$R$
has trivial algebraic intersection with 
any closed loop it contains (since it is orientable).
Once we lift 
$R$
to the cover
$S\times \RR$,
we can use incompressibility and 
$\partial$-incompressibility of
$R$
to finish the proof of the claim.

Using this claim we can prove the proposition. 
The conclusion clearly holds if 
$F$
has at most two summands.
Now
we can assume that every surface
$F_i$
is a fibre of 
$M$
and that there are more than two summands in the 
whole expression. 

Let 
$F_j$
and
$F_k$
be two summands in
$F=k_1F_1+\ldots+k_nF_n$
that have non-trivial intersection. After making
regular alterations along all the curves in
$F_j\cap F_k$
we get
$F_j+F_k=A+B$
where 
$A$
and 
$B$
are disjoint connected normal surfaces that are isotopic to
the fibre of
$M$.
We can see this in the following way. 
Each component 
$D$
of the surface
$F_k+F_j$
appears as a summand in some normal sum representing 
$F$.
It is therefore injective 
$\partial$-incompressible and homologically non-trivial in
$H_2(M,\partial M;\ZZ_2)$
(the surface 
$D$
can be neither separating because 
$\chi(F)<\chi(D)<0$ 
nor can it be closed since it is injective).
Think of
$F_k$
and
$F_j$
as non-trivial elements of 
$H_2(M,\partial M;\ZZ_2)$.
Their normal sum is therefore zero in
$H_2(M,\partial M;\ZZ_2)$
which means that 
$F_j+F_k$
has an even number of components. If there are four or more,
then at most one is a fibre because
$\chi(F_j+F_k)=2\chi(S)$,
where 
$S$
is a fibre of 
$M$.
In that case at least one component is a bounded one-sided surface
with its Euler characteristic strictly larger than
$\frac{1}{2}\chi(S)$.
This is a contradiction since
$\chi(F)<\chi(F_j)+\chi(F_k)= 2\chi(S)$
and therefore the double of that one-sided surface would 
satisfy all the conditions from the proposition 
with its Euler characteristic strictly larger than that of
$F$.
So we can conclude that 
$F_j+F_k=A'\cup B'$,
where  
$A'$
and
$B'$
are disjoint connected homologically non-trivial surfaces.
We also have
$\chi(A')+\chi(B')=2\chi(S)$.
If both 
$A'$
and
$B'$
are one-sided, then the double
of the component with the larger Euler characteristic
implies
$2\chi(S)\leq\chi(F)$,
which is a contradiction. If only one of them
is one-sided, then the other one 
is isotopic to the fibre by the claim.
So the Euler characteristic of 
the double of the one-sided component equals
$2\chi(S)$
which leads to contradiction as before.
So in the end we get that both
components are two-sided and hence fibres by the claim.

If
$F$
is a sum of at least four fibres (there can not be three because 
$F$
is trivial in
$H_2(M,\partial M;\ZZ_2)$) 
then, by what we've just proved, we never
reduce the number of fibres in the sum by doing regular 
alterations. This is a contradiction because 
$F$
is connected. So the proposition follows.
\end{proof}

The only classes of 3-manifolds we need to think about now
in order to finish the proof of theorem~\ref{thm:final} are:
closed atoroidal surface semi-bundles which are not rational homology
3-spheres, closed atoroidal surface bundles with first Betti number at least
2 and an-annular 
semi-bundles with a single boundary component which contain
no closed injective surfaces other than the boundary torus and which are
not homeomorphic to surface bundles. In all these manifolds
we can find a homologically non-trivial surface which is not a fibre
in any fibration of the manifold. It follows
from~\cite{thurston} that any connected incompressible
$\partial$-incompressible surface whose homology
class is carried by a vertex in the boundary of the unit ball 
for the Thurston norm on
$H_2(M,\partial M;\RR)$ 
can not be a fibre. If we pick one such surface which minimises
the topological complexity in its homology class,
then its Euler characteristic is bounded by corollary 5.8
in~\cite{wang}. Since all 3-manifolds on the above list are
both atoroidal and an-annular, we can use this bound and the
techniques developed in subsection 4.2 of~\cite{mijatov2} to control 
the normal complexity of our surface. This completes
the proof of theorem~\ref{thm:final}.

\bibliographystyle{amsplain}
\bibliography{cite}

\noindent \textsc{Department of Pure Mathematics and Mathematical Statistics},
\textsc{Center for Mathematical Sciences},
\textsc{University of Cambridge}\\
\textsc{Wilberforce Road},
\textsc{Cambridge, CB3 0WB},
\textsc{UK}\\
\textit{E-mail address:} \texttt{a.mijatovic@dpmms.cam.ac.uk}

\end{document}